\documentclass[11pt,a4paper,reqno]{amsart}
\usepackage{graphicx}
\usepackage{epsfig}
\usepackage{amsmath}
\usepackage{amsfonts}
\usepackage{amssymb}
\usepackage{amscd}

\begin{document}

	\begin{minipage}[b]{0.5\linewidth}
		{\includegraphics[height=0.82in,width=5.94in]{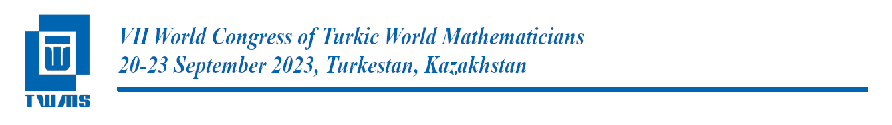}}
	\end{minipage}

\bigskip
\title []
{ON THE MINIMIZATION OF $k$-VALUED LOGIC FUNCTIONS IN THE CLASS OF DISJUNCTIVE NORMAL FORMS}

\author[]{Anvar~Kabulov$^{*~ 1, 3, 4, 5}$, Abdussattar~Baizhumanov $^{2}$, Mansur~Berdimurodov$^{1}$}

\maketitle
\begin{center}
	{\scriptsize $^{1}$National University of Uzbekistan named after Mirzo Ulugbek, Tashkent, Uzbekistan\\
	 \scriptsize $^{2}$South Kazakhstan State Pedagogical University, Shymkent, Kazakhstan\\
	 \scriptsize $^{3}$Tashkent State Technical University, Tashkent, Uzbekistan\\
	 \scriptsize $^{4}$Gulistan State University, Sirdaryo Region, Uzbekistan\\
	 \scriptsize $^{5}$Engineering Federation of Uzbekistan, Tashkent, Uzbekistan\\
	  \indent $^{*}$e-mail: islambeksaymanov@gmail.com }
\end{center}

\begin{abstract}
The paper considers the representation of $k$-valued logical functions in the class of disjunctive normal forms. Various classes of monotone functions of $k$-valued logic are investigated. Theorems are proved on the coincidence of reduced and shortest disjunctive nominal forms of $k$-valued functions. For a certain class of $k$-valued monotone functions, we prove an estimate for the number of functions in this class. we prove criteria for the absorption of elementary conjunctions by a first-order neighborhood of disjunctive normal forms of $k$-valued functions.

\bigskip
\bigskip
{\footnotesize \noindent \textbf{Keywords:} $k$-valued, minimization, disjunctive normal form, rank, abbreviated d.n.f., monotone function.

}

\end{abstract}

\section{Introduction}
The methods of $k$-valued logics are generally necessary for the study of a number of important problems from the most diverse fields: biology, medicine, military affairs, automation, control, planning of experiments etc., everywhere where not only the quantitative relationships between the quantities characterizing the processes under consideration are significant, but also the logical dependences connecting them \cite{bib1,bib2,bib3,bib4}. A multi-valued logical function can be represented as a disjunction (multi-place function "or") ${{K}_{1}}\vee {{K}_{2}}\vee ...\vee {{K}_{m}}$, where each term is a conjunction (multi-place function “and”) of certain variables from the set $\left\{ {{x}_{1}},\text{  }.\text{ }.\text{ }.\text{  },\text{   }{{x}_{n}} \right\}$, taken with or without negation. $k$-valued function gives a description functioning of the control system, and the formula realizing it, in particular, the disjunctive normal form (d.n.f.), describes the scheme of this system, so that the nodes and elements of the scheme correspond to the terms and letters of the d.n.f. as a rule, $k$-valued function has many essentially different d.n.f. \cite{bib5,bib6,bib7,bib8,bib9,bib10,bib11,bib12}. In mathematical logic they are considered from the qualitative side. With the development of cybernetics, the terms and letters of the d.n.f. began to reflect equipment costs in circuits and this drew attention to the quantitative side. Therefore, one of the problems of $k$-valued logics dictated by practice is the problem of minimizing multi-valued functions. The results of research in some areas in this area, in particular, minimization in certain systems multi-valued functions are quite widely displayed in the literature \cite{bib13,bib14,bib15}. Therefore, it should immediately be noted that we will only discuss the minimization of multi-valued functions in the class of d.n.f..

\section{Problem statement} 
Consider the set of multivalued logic functions depending on $n$ variables, etc. the set of functions defined on the set of all vertices of the $n$-dimensional $k$-lattice $E_{n}^{k}$ and taking values from the set ${{E}_{k}}=\left\{ 0,1,\ldots ,k-1 \right\}$. With this interpretation, there is a one-to-one correspondence between multivalued logic functions depending on $n$ arguments and subsets of ${{N}_{f}}\subseteq E_{n}^{k}$. The function $f\left( {{x}_{1}},\ldots ,{{x}_{n}} \right)$ and the subset ${{N}_{f}}\subseteq E_{n}^{k}$ correspond to each other in the case
$f\left( x \right) = \left\{ \begin{array}{l}
\gamma ,\,\,\,if\,\,x\in {{N}_{f}} \\ 
0,\,\,\,if\,\,x\in E_{n}^{k}\backslash {{N}_{f}} \\ 
\end{array} \right.,$ where
\begin{equation}
\gamma \in {{E}_{f}}\subseteq \left\{ 0,1,\ldots ,k-1 \right\}
\end{equation}

We can assume that the set ${{E}_{f}}$ divides the function $f$ into a number of subfunctions ${{f}_{{{\gamma }_{1}}}}\left( \widetilde{x} \right),\ldots ,{{f}_{{{\gamma }_{m}}}}\left( \widetilde{x} \right)$, and the set ${{N}_{f}}$ into pairwise disjoint subsets ${{N}_{{{f}_{{{\gamma }_{1}}}}}},\ldots ,{{N}_{{{f}_{{{\gamma }_{m}}}}}}$, where $m=\left| {{E}_{f}} \right|$,
$${{f}_{{{\gamma }_{i}}}}\left( \widetilde{x} \right)=\begin{cases}
{{\gamma }_{i}},\,\,if\,\,f\left( x \right)={{\gamma }_{i}} \\ 
0,\,\,if\,\,f\left( x \right)\ne {{\gamma }_{i}} \\ 
\end{cases} $$
$${{N}_{{{f}_{i}}}}=\left\{ \widetilde{\alpha }:\left( \widetilde{\alpha }\in E_{n}^{k} \right)\wedge \left( f\left( \widetilde{\alpha } \right)={{\gamma }_{i}} \right),\left( i=\overline{1,m} \right). \right\}$$
~~~It is easy to see that
\begin{equation}
f\left( \widetilde{\alpha } \right)=\underset{{{\gamma }_{i}}}{\mathop{\max }}\,\left\{ {{f}_{{{\gamma }_{1}}}}\left( \widetilde{\alpha } \right),\ldots ,{{f}_{{{\gamma }_{m}}}}\left( \widetilde{\alpha } \right) \right\}
\end{equation}

The function ${{f}_{\gamma }}\left( \widetilde{x} \right)$, which takes only two values $\left( 0\,\,and\,\,\gamma  \right)$, will be called quasi-Boolean, and the representation of the function $f\left( {{x}_{1}},\ldots ,{{x}_{n}} \right)$ in the form of a (2)-quasi-Boolean representation of the multi-valued logical function $f\left( {{x}_{1}},\ldots ,{{x}_{n}} \right)$.

We introduce the function

${{J}_{M}}( x )= \begin{cases}
k-1,\,\,\,if\,\,x\in M \\ 
0,\,\,\,\,\,\,\,\,\,\,if\,\,x\notin M \\ 
\end{cases} $, where
\begin{equation}
M\subseteq {{E}_{k}}=\left\{ 0,1,2,\ldots ,k-1 \right\}
\end{equation}
An elementary conjunction (e.c.) is an expression

$\mathfrak{A}=\min \left[ {{J}_{{{M}_{1}}}}\left( {{x}_{1}} \right),\ldots ,{{J}_{{{M}_{n}}}}\left( {{x}_{n}} \right),\gamma  \right]$, where \begin{equation}
\varnothing \ne {{M}_{j}}\subseteq {{E}_{k}},\left( j=\overline{1,n} \right)
\end{equation}
Further, for brevity, formulas $\max \left[ {{\mathfrak{A}}_{1}},\ldots ,{{\mathfrak{A}}_{m}} \right]$ will be conventionally denoted as ${{\mathfrak{A}}_{1}}\vee \ldots \vee {{\mathfrak{A}}_{m}}=\underset{i=1}{\overset{m}{\mathop{\vee }}}\,{{\mathfrak{A}}_{i}}:$ if ${{\mathfrak{A}}_{i}}$ is an analog of e.c., then this formula will be called disjunctive normal form (d.n.f.).

The area of truth of e.c. let's call $\mathfrak{A}$ the region ${{N}_{\mathfrak{A}}}$ in which $\mathfrak{A}$ takes the value $\gamma $. It is easy to see that the domain ${{N}_{\mathfrak{A}}}=\prod\limits_{j=1}^{n}{{{M}_{j}}}$ is a sub lattice (a subset of the set $E_{n}^{k}$) of the lattice $E_{n}^{k}$. With such a geometric consideration, the e.c. the sub lattice ${{N}_{\mathfrak{A}}}$ corresponds to the lattice $E_{n}^{k}$.

Rank e.c. let's $\mathfrak{A}$ call a number $r\left( \mathfrak{A} \right)=\sum\limits_{j=1}^{n}{\left( k-\left| {{M}_{j}} \right| \right)}=kn-\sum\limits_{j=1}^{n}{\left| {{M}_{j}} \right|}$.

The formula $\mathfrak{M}=\underset{i=1}{\overset{t}{\mathop{\vee }}}\,{{\mathfrak{A}}_{i}}$ where all ${{\mathfrak{A}}_{i}},\,\left( i=\overline{1,t} \right)$ are e.c. will be called the disjunctive normal form.

Note that each set-valued logic function $f\left( {{x}_{1}},\ldots ,{{x}_{n}} \right)$ corresponds to a non-empty class of d.n.f. realizing the given function. The set of all intervals corresponding to e.c. a certain d.n.f. from this class determines the covering of ${{N}_{f}}$ by sub lattices of the lattice $E_{n}^{k}$. Hence it follows that the subsets $M\subseteq E_{n}^{k}$ can be defined using the d.n.f.

Let $I=\left\{ {{N}_{\mathfrak{A}}} \right\}$, be some subset of sub lattices from $E_{n}^{k}$. 

A sub lattice ${{N}_{B}}\in I$ is said to be maximal with respect to $M$ if there is no sub lattice ${{N}_{\mathfrak{A}}}$ in $I$ such that ${{N}_{\mathfrak{A}}}\ne {{N}_{B}}$ and ${{N}_{\mathfrak{A}}}\supseteq {{N}_{B}}$.

To represent the function $f\left( {{x}_{1}},\ldots ,{{x}_{n}} \right)$ in the form of a d.n.f. we considered its quasi-Boolean representation: $f={{f}_{{{\gamma }_{1}}}}\vee \ldots \vee {{f}_{{{\gamma }_{m}}}}$ and ${{\gamma }_{1}}<{{\gamma }_{2}}<\ldots <{{\gamma }_{m}}$.

Note that for the same function $f\left( \widetilde{x} \right)$ there can be several equivalent quasi-Boolean representations. Indeed, $f={{f}_{{{\gamma }_{1}}}}\vee \ldots \vee {{f}_{{{\gamma }_{m}}}}={{f}^{*}}=f_{{{\gamma }_{1}}}^{*}\vee \ldots \vee f_{{{\gamma }_{m}}}^{*}$ where ${{N}_{f_{{{\gamma }_{i}}}^{*}}}={{N}_{{{f}_{{{\gamma }_{i}}}}}}\cup {{Q}_{i}}$, ${{Q}_{i}}\subseteq \underset{j>i}{\mathop{\cup }}\,{{N}_{{{\gamma }_{j}}}},\left( i=\overline{1,m} \right)$.

We will consider only one "maximum" representation ${f}'={{{f}'}_{{{\gamma }_{1}}}}\vee \ldots \vee {{{f}'}_{{{\gamma }_{m}}}}$ where ${{N}_{{{{{f}'}}_{{{\gamma }_{i}}}}}}=\underset{j=1}{\overset{n}{\mathop{\cup }}}\,{{N}_{{{f}_{{{\gamma }_{j}}}}}},\,\left( i=\overline{1,m} \right)$.

Select all maximal sub lattices ${{N}_{B_{j}^{i}}},\left( i=\overline{1,m} \right)$ contained in ${{N}_{{{{{f}'}}_{{{\gamma }_{i}}}}}}$ that have non-empty intersection with ${{N}_{{{f}_{{{\gamma }_{i}}}}}}$ and such that the value of $B_{j}^{i}$ is equal to ${{\gamma }_{i}}$ in ${{N}_{B_{j}^{i}}},\left( i=\overline{1,m} \right)$. D.n.f. $\mathfrak{M}=\underset{i=1}{\overset{m}{\mathop{\vee }}}\,\underset{j=1}{\overset{{{I}_{i}}}{\mathop{\vee }}}\,B_{j}^{i}$ is called the reduced disjunctive normal form of the function $f\left( \widetilde{x} \right)$.

A covering of a set ${{N}_{f}}$ by maximal sub lattices is said to be irreducible if, after the removal of any of its sub lattices, it ceases to be a covering. A d.n.f. realizing a function $f$ is called dead-end if it corresponds to an irreducible cover of the set ${{N}_{f}}$.

Consider a multivalued logic function $F\left( {{x}_{1}},\ldots ,{{x}_{n}} \right)$ defined at $M\subseteq E_{n}^{k}:F\left( \widetilde{x} \right)={{\gamma }_{j}}$, if $\widetilde{x}\in {{M}_{j}},\left( j=\overline{1,m} \right)$, $m<k,\,{{\gamma }_{j}}\in {{E}_{k}},\,M=\underset{i=0}{\overset{m}{\mathop{\bigcup }}}\,{{M}_{i}}$ and ${{M}_{i}}\bigcap {{M}_{j}}=\varnothing $ at $\left( i\ne j,\,i,\,j=\overline{0,m} \right)$. And ${{\gamma }_{1}}<\ldots <{{\gamma }_{m}},\,{{\gamma }_{0}}=0$.

Thus, $F\left( {{x}_{1}},\ldots ,{{x}_{n}} \right)$ is defined by specifying pairwise disjoint sets ${{M}_{0}},\ldots ,{{M}_{m}}$. The function $F\left( \widetilde{x} \right)$ is defined, generally speaking, not on the entire set $E_{n}^{k}$. There are different before the definition in the class of functions $F\left( \widetilde{x} \right)$, multi-valued logic, not equivalent to each other.

Our task is to find the simplest ones, in a certain sense, before definitions.

For $F\left( \widetilde{x} \right)$, select all maximal intervals ${{N}_{B_{j}^{i}}}$, $\left( i=\overline{1,m},\,j=\overline{0,{{I}_{i}}} \right)$ contained in $E_{n}^{k}\backslash \bigcup\limits_{v=0}^{i-1}{{{M}_{v}}}$ that have non-empty intersection with ${{M}_{i}}$ such that the value of $B_{j}^{i}$ is equal to ${{\gamma }_{i}}$.

D.n.f. $\mathfrak{M}=\underset{i=1}{\overset{m}{\mathop \vee }}\,\underset{j=1}{\overset{{{I}_{i}}}{\mathop \vee }}\,\,B_{j}^{i}$ is called the reduced normal form for $F\left( \widetilde{x} \right)$. It is easy to see that d.n.f. ${{\eta }_{\Sigma TF}}$ is uniquely determined by the function $F$.

Let us now indicate the points at which, when the values of the function $F$ change (transition to ${F}'$), the values of ${{\eta }_{\Sigma TF}}$ change (transition to ${{\eta }_{\Sigma T{F}'}}$).

\section{Monotone functions of $k$-valued logic} 

Let's consider some order on the set ${{\varepsilon }_{k}}$. For two sets $\Im =\left( {{\alpha }_{1}},{{\alpha }_{2}},\ldots {{\alpha }_{n}} \right)$ and $\widetilde{\beta }=\left( {{\beta }_{1}},{{\beta }_{2}},\ldots {{\beta }_{n}} \right)$ the precedence relation $\Im \le \widetilde{\beta }$ is satisfied if ${{\alpha }_{i}}\le {{\beta }_{i}}$ in this order for any $i=\overline{1,n}$.

\textbf{Definition~1.} A $k$-valued logic function $f\left( {{x}_{1}},{{x}_{2}},\ldots ,{{x}_{n}} \right)$ is said to be monotonic with respect to a given order if for any tuples $\alpha $ and $\beta $ such that $\left( \Im \le \widetilde{\beta } \right)$ we have $f\left( \alpha  \right)\le f\left( \beta  \right)$.

If $0<1<2<\ldots <k-1$, then the set of functions that are monotone in this order constitutes the class of monotone functions of $k$-valued logic.

\textbf{Theorem~1.} Abbreviated d.n.f. monotone function ${{f}_{\gamma }}$ of $k$-valued logic in $n$ variables

a)	consists of e.c. ${{K}^{A}}$, and only elementary formulas of the form ${{J}_{\left[ a,k-1 \right]}}\left( x \right),0\le a\le k-1$ are used;

b)	Is the only minimal (shortest) d.n.f. of the function $f$.

\textbf{Proof.} 

a)~Let $K={{J}_{{{T}_{1}}}}\left( {{x}_{1}} \right)\cdot {{J}_{{{T}_{2}}}}\left( {{x}_{2}} \right)\cdot \ldots \cdot {{J}_{{{T}_{n}}}}\left( {{x}_{n}} \right)\cdot \gamma $, where ${{T}_{j}}=\left[ {{a}_{j}},k-1 \right]$, $j\ne i,\,0\le {{a}_{j}}\le k-1,\,\,{{T}_{i}}=\left[ {{b}_{i}} \right]\bigcup \left[ {{a}_{i}},k-1 \right],\,{{b}_{i}}\le {{a}_{i}}$.

Then the conjunction $K$, and hence the function ${{f}_{\gamma }}$, takes the value $\gamma $ on the set $\widetilde{a}=\left( {{a}_{1}},{{a}_{2}},\ldots ,{{a}_{i-1}},{{a}_{i}},{{a}_{i+1}},\ldots ,{{a}_{n}} \right)$.

b)~It follows from the monotonicity condition for the function ${{f}_{\gamma }}$ that for any set $\widetilde{b}$ such that $\widetilde{b}\ge \widetilde{a},\,{{f}_{\gamma }}\left( \widetilde{b} \right)=\gamma $ therefore, there exists an e.c. ${K}'={{J}_{{{{{T}'}}_{1}}}}\left( {{x}_{1}} \right)\cdot {{J}_{{{{{T}'}}_{2}}}}\left( {{x}_{2}} \right)\cdot \ldots \cdot {{J}_{{{{{T}'}}_{n}}}}\left( {{x}_{n}} \right)\cdot \gamma $, where ${{{T}'}_{j}}={{T}_{j}},\,j\ne i,\,{{{T}'}_{i}}=\left[ {{b}_{i}},k-1 \right]$, for which ${{U}_{K}}\subset {{U}_{{{K}'}}}\subseteq {{U}_{{{f}_{\gamma }}}}$, and e.c. $K$ is not maximal for ${{U}_{{{f}_{\gamma }}}}$.

c)~As proved, each maximum e.c. $K$ functions has the form $K={{J}_{\left[ {{a}_{1}},k-1 \right]}}\left( {{x}_{1}} \right)\cdot {{J}_{\left[ {{a}_{2}},k-1 \right]}}\left( {{x}_{2}} \right)\cdot \ldots \cdot {{J}_{\left[ {{a}_{n}},k-1 \right]}}\left( {{x}_{n}} \right)\cdot \gamma $,  $0\le {{a}_{j}}\le k-1,\,\,j=\overline{1,n}$.

Let us show that the set $\widetilde{a}=\left( {{a}_{1}},{{a}_{2}},\ldots ,{{a}_{n}} \right)$ is core for the function ${{f}_{\gamma }}$, etc. in the abbreviated d.n.f. the function ${{f}_{\gamma }}$ has no e.c., except for $K$, which takes the value $\gamma $ on this set.

Indeed, if in the reduced d.n.f. function ${{f}_{\gamma }}$ was an e.c. ${K}'$, which takes the values $\gamma $ on the tuple $\widetilde{a}$, then it would follow from the monotonicity condition for the function ${{f}_{\gamma }}$ that the e.c. ${K}'$ takes the value $\gamma $ on all tuples $\widetilde{b}$ such that $\widetilde{b}\ge \widetilde{a}$, then ${{U}_{K}}\subseteq {{U}_{{{K}'}}}$, which contradicts the maximum e.c. $K$. 

\textbf{The theorem is proved.}

\textbf{Corollary.} Abbreviated d.n.f. monotone function $f$ of $k$-valued logic in n variables consists of e.c. ${{K}^{A}}$ and is the only minimal d.n.f. functions $f$.

\textbf{Proof.} A monotone function $f$ has the following obvious properties: for any comparable collections $\widetilde{a}\in {{N}_{{{f}_{\gamma }}}}$ and $\widetilde{b}\in {{N}_{{{f}_{\theta }}}}$ we have $\widetilde{b}>\widetilde{a}$ for $\gamma <\theta $. Therefore, e.c. $K$, is included in the abbreviated d.n.f. functions ${{f}_{\gamma }}\left( \gamma \in \left\{ {{\varepsilon }_{k}}/0 \right\} \right)$ do not contain elementary formulas of type ${{J}_{T}}\left( x \right)$, where the set $T$ is a disconnected set of points from ${{\varepsilon }_{k}}$. 

The second assertion is proved by the method of theorem 1, and the core sets for the function $f$ are the sets $\widetilde{a}=\left( {{a}_{1}},{{a}_{2}},\ldots ,{{a}_{n}} \right)$ in the e.c. $K={{J}_{\left[ {{a}_{1}},{{b}_{1}} \right]}}\left( {{x}_{1}} \right)\cdot {{J}_{\left[ {{a}_{2}},{{b}_{2}} \right]}}\left( {{x}_{2}} \right)\cdot \ldots \cdot {{J}_{\left[ {{a}_{n}},{{b}_{n}} \right]}}\left( {{x}_{n}} \right)\cdot \gamma $. 

\textbf{Corollary proven.} On the set ${{\varepsilon }_{k}}$ we introduce a partial order $0<1,\,0<2,\,\ldots ,0<k-1$ where $i$ is incomparable with $j$ if $i,j\in \left\{ {{\varepsilon }_{k}}/0 \right\}$.

The set of functions of $k$-valued logic $f$ in $n$ variables, monotone in a given order, is combined into the class $S$. Let us estimate the cardinality of the class $S$.

The set ${{\varepsilon }_{k}}$ is associated with a basic graph - a directed graph with $K$ vertices corresponding to the elements of the set ${{\varepsilon }_{k}}$, in which there is an arc $\left( i,j \right)$ if and only if $i>j$.

Let us introduce the $Z$ axis on the plane. Associate each point $A$ of the plane with the number ${{Z}_{A}}$, the projection of $A$ onto the $Z$ axis. In particular, if the basis graph is drawn on the plane, then each vertex corresponds to the number ${{Z}_{i}}$. An image of a basic graph is called admissible if for any arc $\left( i,j \right)$ of the graph ${{Z}_{i}}-{{Z}_{j}}\ge 1$.

Consider a random variable $\xi ={{Z}_{A}}$, where the point $A$ can fall into any vertex of the graph with probability $\frac{1}{k}$ then the expectation is $M \xi ={{Z}_{cp}}=\frac{{{Z}_{1}}+\ldots +{{Z}_{k}}}{k}$ and the variance is $D\xi =\frac{1}{k}\sum\limits_{i=1}^{k}{{{\left( {{Z}_{i}}-{{Z}_{cp}} \right)}^{2}}}$.

We will consider the image of the graph shifted so that $M\xi =0$.

In older articles an estimate is obtained for finding the number $\psi \left( n \right)$ of monotone functions of $n$ variables from an arbitrary partially ordered set of $k$ elements:
\begin{equation}
\psi \left( n \right)={{d}^{\frac{1}{\sqrt{2\pi D}}\cdot \frac{{{k}^{n}}}{\sqrt{n}}\left( 1+\varepsilon \left( n \right) \right)}}
\end{equation}
where $\varepsilon \left( n \right)\to 0$ at $n\to \infty $; $D=\inf \,D\xi $; $d=\max \left( \left| {{H}_{1}} \right|,\ldots ,\left| {{H}_{s}} \right| \right)$, ${{H}_{0}}\le {{H}_{1}}\le \ldots \le {{H}_{s+1}}$, all ${{H}_{i}}\subseteq {{\varepsilon }_{k}},\,s\ge 1$,  $\left| {{H}_{0}} \right|=\left| {{H}_{s+1}} \right|=1$ and ${{H}_{i}}\ne {{H}_{j}}$, at $i\ne j$ (${{H}_{i}}\le {{H}_{j}}$ if $a\le b$ for any $a\in {{H}_{i}},\,\,b\in {{H}_{j}}$), the maximum is taken over all possible chains. This estimate is also valid for the class of functions $S$.

For the order in ${{\varepsilon }_{k}}$ we have:
\begin{equation} \begin{cases}
{{Z}_{1}}-{{Z}_{0}}\ge 1 \\ 
\cdots \,\,\,\cdots \,\,\,\cdots  \\ 
{{Z}_{k-1}}-{{Z}_{0}}\ge 1 \\ 
\end{cases} \end {equation}
\begin{equation}
{{Z}_{cp}}=\frac{{{Z}_{0}}+{{Z}_{1}}+\ldots +{{Z}_{k-1}}}{k}=0,\,D\xi =\frac{1}{k}\sum\limits_{i=0}^{k-1}{{{Z}_{i}}^{2}}
\end{equation}
It follows that ${{Z}_{0}}\le -\frac{k-1}{k}$ and ${{Z}_{i}}\ge \frac{1}{k}$, $i=\overline{1,k-1}$, so $D\xi \ge \frac{k-1}{{{k}^{2}}}$ and $D\ge \frac{k-1}{{{k}^{2}}}$. 

On the set ${{\varepsilon }_{k}}$ for the introduced order there is only $\left( k-1 \right)$ chain
\begin{equation}
\begin{cases}
\left\{ 0 \right\}<\left\{ 0,1 \right\}<\left\{ 1 \right\} \\ 
\cdots \cdots \cdots  \\ 
\left\{ 0 \right\}<\left\{ 0,k-1 \right\}<\left\{ k-1 \right\} \\ 
\end{cases}
\end{equation}

It is obvious that in the considered case $d=2$.

Consequently 
\begin{equation}
\begin{cases}
\psi \left( n \right)={{d}^{\frac{1}{\sqrt{2\pi \left( k-1 \right)}}\cdot \frac{{{k}^{n+1}}}{\sqrt{n}}\left( 1+\varepsilon \left( n \right) \right)}}
\end{cases}
\end{equation}

where $\varepsilon \left( n \right)\to 0$ at $n\to \infty $.

For functions of class $S$, e.c. consists of elementary formulas of the form ${{J}_{T}}\left( x \right)$, where $T\subseteq \left\{ {{\varepsilon }_{k}}/0 \right\}$.

If for monotone functions $f$ of $k$-valued logic the abbreviated d.n.f. is the only minimal one, then for functions of the class $S$ this property does not hold, which is demonstrated by the following.

\textbf{Example.} $k=3,\,n=3.$

Let the function $f\left( {{x}_{1}},{{x}_{2}},{{x}_{3}} \right)$ take the values 1 on the sets $\left( 0,1,1 \right),\left( 1,1,1 \right),\left( 1,2,1 \right),\left( 2,1,1 \right),\left( 1,2,2 \right)$ and 0 in other cases.

Abbreviated d.n.f. function $f$ has the form:
\begin{equation}
{{D}_{C}}\left( f \right)={{J}_{1}}\left( {{x}_{2}} \right)\cdot {{J}_{1}}\left( {{x}_{3}} \right)\vee {{J}_{1}}\left( {{x}_{1}} \right)\cdot {{J}_{2}}\left( {{x}_{2}} \right)\cdot {{J}_{[1,2]}}\left( {{x}_{3}} \right)
\end{equation}
and the minimum:
\begin{equation}
{{D}_{M}}\left( f \right)={{J}_{1}}\left( {{x}_{2}} \right)\cdot {{J}_{1}}\left( {{x}_{3}} \right)\vee {{J}_{1}}\left( {{x}_{1}} \right)\cdot {{J}_{[1,2]}}\left( {{x}_{2}} \right)\cdot {{J}_{1}}\left( {{x}_{3}} \right)
\end{equation}
The process of transition from the abbreviated d.n.f. functions $f$ of a $k$-valued logic to a dead-end one can be divided into elementary steps, each of which is a removal from the d.n.f. $D$ obtained in the previous step, one e.c. $K$. Removed e.c. is such that ${{U}_{K}}\subseteq \bigcup\limits_{j=1}^{m}{{{U}_{{{K}_{j}}}}}$, where ${{K}_{j}}$ are some e.c. from d.s.f. $D$ different from $K$.

In older articles the criterion for covering an interval by the sum of other intervals for functions $f$ of $k$-valued logic is described. For functions $f$ of class $S$, this criterion has a simpler form.

E.c. ${{K}_{1}}$ and ${{K}_{2}}$ are called orthogonal if ${{K}_{1}}\cdot {{K}_{2}}=0$. In other words, conjunctions ${{K}_{1}}$ and ${{K}_{2}}$ are orthogonal if and only if ${{U}_{{{K}_{1}}}}\bigcap {{U}_{{{K}_{2}}}}=\varnothing $. Obviously, when studying the absorption, some sets of e.c. $\left\{ {{K}_{j}} \right\},\,j=\overline{1,m}$ e.c. $K$ it suffices to consider only those e.c. that are non-orthogonal to $K$. To check orthogonality, it is easiest to use the following properties: two e.c. are orthogonal if and only if there exists a variable ${{x}_{j}}$ for which $T_{j}^{1}\bigcap T_{j}^{2}=\varnothing $ is satisfied in the elementary formulas ${{J}_{T_{j}^{1}}}\left( {{x}_{1}} \right)$ and ${{J}_{T_{j}^{2}}}\left( {{x}_{2}} \right)$.

D.n.f. $D$ realizing the function f absorbs the e.c. $K$ if $K\left( \widetilde{x} \right)\le D\left( \widetilde{x} \right)$ for any $\widetilde{x}=\varepsilon _{k}^{n}$. 

So let $K={{J}_{{{T}_{1}}}}\left( {{x}_{1}} \right)\cdot {{J}_{{{T}_{2}}}}\left( {{x}_{2}} \right)\cdot \ldots \cdot {{J}_{{{T}_{t}}}}\left( {{x}_{t}} \right)\cdot \gamma $.

Obviously, $K$ can be absorbed only by those sets of e.c. $\left\{ {{K}_{j}} \right\}$, which take values from $\left\{ 0,\gamma  \right\}$, so let's consider the absorption process using the example of the quasi-Boolean function ${{f}_{\gamma }}$.

For each e.c. $\left\{ {{K}_{j}} \right\},\,\,j=\overline{1,m}$ construct an e.c. ${{K}_{j}}$, replacing the elementary formulas ${{J}_{T}}\left( x \right)$ occurring in ${{K}_{j}}$ with ${{J}_{[1,k-1]}}\left( x \right)$.

It's obvious that ${{U}_{{{K}_{j}}}}\subseteq {{U}_{{{K}_{j}}}},\,\,{{U}_{{{K}_{j}}}}\cap \left( {{U}_{D}}/{{U}_{{{K}_{j}}}} \right)=\varnothing $.

Let us introduce into consideration the set $\varepsilon _{k}^{n,t}$- the collection of all sets from $\varepsilon _{k}^{n}$, in which the $t$ first coordinates take values from $\left\{ {{\varepsilon }_{k}}/0 \right\}$, and the rest are arbitrary.

\textbf{Theorem~2.} The disjunction $D=\underset{j=1}{\overset{m}{\mathop{\vee }}}\,{{K}_{j}}$ absorbs the e.c. $K$ if and only if $\underset{j=1}{\overset{m}{\mathop{\vee }}}\,{{K}_{j}}=\gamma $ for any $\widetilde{x}=\varepsilon _{k}^{n,t}$, etc. if $\underset{j=1}{\overset{m}{\mathop{\vee }}}\,{{K}_{j}}=\varepsilon _{k}^{n,t}$.

\textbf{Proof.} \textit{Need.} Let $D$ absorb the e.c. $K$. Let us prove that in this case $\underset{j=1}{\overset{m}{\mathop{\vee }}}\,{{K}_{j}}=\gamma $ for any $\widetilde{x}=\varepsilon _{k}^{n,t}$. Let us assume that this is not the case, etc. there is a collection $\Im $ such that $\underset{j=1}{\overset{m}{\mathop{\vee }}}\,{{K}_{j}}\left( \Im  \right)=0$. Denote by ${{x}_{{{i}_{1}}}},\ldots ,{{x}_{{{i}_{p}}}}$ the variables that are not included in any of the e.c. from $D$. Obviously, the values of the remaining variables do not affect the value of the expression $\underset{j=1}{\overset{m}{\mathop{\vee }}}\,{{K}_{j}}$.

The value of the function $f$ on the tuples $\Im $ will be denoted by $\left[ f \right]\Im $. Then one can write
\begin{equation}
\left[ \underset{j=1}{\overset{m}{\mathop{\vee }}}\,{{K}_{j}} \right]\left\{ \Im  \right\}=0
\end{equation}
where $\left\{ \Im  \right\}$ is the set of sets whose $\left( {{x}_{1}},{{x}_{2}},\ldots ,{{x}_{t}} \right)$ coordinates take all possible values from $\left\{ {{\varepsilon }_{k}}/0 \right\}$, and the remaining $\left( n-t \right)$ coordinates are such that (12) is satisfied. hence we get that $\left[ {{K}_{j}} \right]\left\{ \Im  \right\}=0$ for all $j=\overline{1,m}$, so $\left[ D \right]\Im =0$.

Let us determine the values of the remaining variables entering $K$ ($K\ne 0$ on ${{U}_{D}}$, since $K$ is not orthogonal to $D$) so that $K$ turns into $\gamma $ on these sets. The intersection of these two sets determines the values of all variables in such a way that $D$ takes the value 0 on this set, and e.c. $K$ value $\gamma $. This contradicts the condition $K\left( \widetilde{x} \right)\le D\left( \widetilde{x} \right)$ for any $\widetilde{x}=\varepsilon _{k}^{n}$, hence the assumption that $\underset{j=1}{\overset{m}{\mathop{\vee }}}\,{{K}_{j}}\ne \gamma $ on $\varepsilon _{k}^{n,t}$ is false, and the necessity of the condition of the theorem is proved.

\textbf{Adequacy.} Let the condition of the theorem be satisfied etc., $\bigcup\limits_{j=1}^{m}{{{U}_{{{K}_{j}}}}}=\varepsilon _{k}^{n,t}$. Then ${{U}_{K}}\subseteq \bigcup\limits_{j=1}^{m}{{{U}_{{{K}_{j}}}}}$, but ${{U}_{K}}\bigcap \left( \bigcup\limits_{j=1}^{m}{{{U}_{{{K}_{j}}}}/\bigcup\limits_{j=1}^{m}{{{U}_{{{K}_{j}}}}}} \right)=\varnothing $. Consequently. ${{U}_{K}}\subseteq \bigcup\limits_{j=1}^{m}{{{U}_{{{K}_{j}}}}}$. 

\textbf{The theorem has been proven.}

\section{Conclusions}
The paper proposes a representation of $k$-valued functions in the class of disjunctive normal forms. Monotone functions of $k$-valued logic are investigated. We prove theorems on the coincidence of abbreviated and shortest d.n.f. $k$-valued functions. For a certain class of $k$-valued monotone functions, the number of functions from this class is calculated. Criteria for the absorption of elementary conjunctions by a first-order neighborhood are proposed.

\end{document}